\documentclass[11pt,reqno]{amsart}
\usepackage{amsmath,amssymb}
\usepackage[colorlinks=true,linkcolor=magenta,citecolor=blue]{hyperref}
\usepackage[margin=1in]{geometry}
\usepackage{enumitem} 
\usepackage{mathabx}
\usepackage{epsfig}

\usepackage{graphicx}
\usepackage{mathrsfs}
\usepackage{color}
\usepackage{amsfonts}
\usepackage{amsmath,amsxtra,latexsym,amsthm, amssymb, amscd, pb-diagram}

\newtheorem{dn}{Definition}[section]
\newtheorem{dl}{Theorem}[section]
\newtheorem{md}{Proposition}[section]
\newtheorem{bd}{Lemma}[section]
\newtheorem{hq}{Corollary}[section]
\newtheorem{nx}{Remark}[section]
\newtheorem{vd}{Example}[section]
\newcommand{\R}{\mathbb{R}}

\newcommand{\Z}{\mathbb{Z}}

\newcommand{\e}{\varepsilon}
\newcommand{\ity}{\infty}

\newcommand{\bbd}{\begin{bd}}
\newcommand{\ebd}{\end{bd}}
\newcommand{\bdn}{\begin{dn}}
\newcommand{\edn}{\end{dn}}
\newcommand{\bhq}{\begin{hq}}
\newcommand{\ehq}{\end{hq}}
\newcommand{\bdl}{\begin{dl}}
\newcommand{\edl}{\end{dl}}
\newcommand{\bnx}{\begin{nx}}
\newcommand{\enx}{\end{nx}}
\newcommand{\bmd}{\begin{md}}
\newcommand{\emd}{\end{md}}
\newcommand{\bvd}{\begin{vd}}
\newcommand{\evd}{\end{vd}}

\title[A blow-up result for semi-linear structurally damped $\sigma$-evolution equations]{A blow-up result for semi-linear structurally damped $\sigma$-evolution equations}
\author{Tuan Anh Dao}
\address{Tuan Anh Dao \hfill\break
$\quad$ School of Applied Mathematics and Informatics, Hanoi University of Science and Technology, No.1 Dai Co Viet road, Hanoi, Vietnam \hfill\break
Faculty for Mathematics and Computer Science, TU Bergakademie Freiberg, Pr\"{u}ferstr. 9, 09596, Freiberg, Germany}
\email{anh.daotuan@hust.edu.vn}
\author{Michael Reissig}
\address{Michael Reissig \hfill\break
$\quad$ Faculty for Mathematics and Computer Science, TU Bergakademie Freiberg, Pr\"{u}ferstr. 9, 09596, Freiberg, Germany}
\email{reissig@math.tu-freiberg.de}

\thanks{The PhD study of MSc. T.A. Dao is supported by Vietnamese Government's Scholarship (Grant number: 2015/911).}

\begin{document}
\subjclass[2010]{35B33, 35B44, 35L76}
\keywords{$\sigma$-evolution equations; structural damping; critical exponent; blow-up; test functions}
	
\begin{abstract}
We would like to prove a blow-up result for semi-linear structurally damped $\sigma$-evolution equations, where $\sigma \ge 1$ and $\delta\in [0,\sigma)$ are assumed to be any fractional numbers. To deal with the fractional Laplacian operators $(-\Delta)^\sigma$ and $(-\Delta)^\delta$ as well-known non-local operators, in general, it seems difficult to apply the standard test function method directly. For this reason, in this paper we shall construct new test functions to overcome this difficulty.
\end{abstract}

\maketitle

\section{Introduction} \label{Sec.main}
The main goal of this paper is to discuss the critical exponent to the following Cauchy problem for semi-linear structurally damped $\sigma$-evolution models:
\begin{equation} \label{pt9.1}
\begin{cases}
u_{tt}+ (-\Delta)^\sigma u+ (-\Delta)^{\delta} u_t=|u|^p, \\
u(0,x)= u_0(x),\,\,\, u_t(0,x)=u_1(x),
\end{cases}
\end{equation}
with some $\sigma \ge 1$, $\delta\in [0,\sigma)$ and a given real number $p>1$. Here, critical exponent $p_{crit}=p_{crit}(n)$ means that for some range of admissible $p>p_{crit}$ there exists a global (in time) Sobolev solution for small initial data from a suitable function space. Moreover, one can find suitable small data such that there exists no global (in time) Sobolev solution if  $1< p \le p_{crit}$. In other words, we have, in general, only local (in time) Sobolev solutions under this assumption for the exponent $p$.
\par For the local existence of Sobolev solutions to (\ref{pt9.1}), we address the interested readers to Proposition $9.1$ in the paper \cite{DabbiccoEbert}. The proof of blow-up results in the present paper is based on a contradiction argument by using the test function method. The test function method is not influenced by higher regularity of the data. For this reason, we restrict ourselves to the critical exponent to (\ref{pt9.1}) in the case, where the data are supposed to belong to the energy space.
In this paper, we use the following notations.
\begin{itemize}[leftmargin=*]
\item For given nonnegative $f$ and $g$ we write $f\lesssim g$ if there exists a constant $C>0$ such that $f\le Cg$. We write $f \approx g$ if $g\lesssim f\lesssim g$.
\item As usual, $H^a$ with $a \ge 0$ stands for Bessel potential spaces based on $L^2$.
\item We denote by $[b]$ the integer part of $b \in \R$.
\item Moreover, we introduce the following two parameters:
$$\mathtt{k}^-:= \min\{\sigma;2\delta\}\quad \text{ and }\quad \mathtt{k}^+:= \max\{\sigma;2\delta\} \quad \text{ if }\delta \in [0,\sigma). $$
\end{itemize}
\par In order to state our main result, we recall the global (in time) existence result of small data energy solutions to (\ref{pt9.1}) in the following theorem.
\bdl[\textbf{Global existence}] \label{dlexistence9.1}
Let $m \in [1,2)$ and $n> m_0 \mathtt{k}^-$ with $\frac{1}{m_0}=\frac{1}{m}- \frac{1}{2}$. We assume the conditions
\begin{align*}
&\frac{2}{m} \le p < \ity &\qquad \text{ if }&\, n \le 2\mathtt{k}^+, \\
&\frac{2}{m} \le p \le \frac{n}{n- 2\mathtt{k}^+} &\qquad \text{ if }&\, n \in \Big(2\mathtt{k}^+, \frac{4\mathtt{k}^+}{2-m}\Big].
\end{align*}
Moreover, we suppose the following condition:
\begin{equation}
p> 1+\frac{m(\mathtt{k}^+ +\sigma)}{n- m\mathtt{k}^-}. \label{exponent9.1}
\end{equation}
Then, there exists a constant $\e_0>0$ such that for any small data
$$ (u_0,u_1) \in \big(L^m \cap H^{\mathtt{k}^+}\big) \times \big(L^m \cap L^2\big) \text{ satisfying the assumption }\|u_0\|_{L^m \cap H^{\mathtt{k}^+}}+ \|u_1\|_{L^m \cap L^2} \le \e_0, $$
we have a uniquely determined global (in time) small data energy solution
$$ u \in C([0,\ity),H^{\mathtt{k}^+})\cap C^1([0,\ity),L^2) $$
to (\ref{pt9.1}). Moreover, the following estimates hold:
\begin{align*}
\|u(t,\cdot)\|_{L^2}& \lesssim (1+t)^{-\frac{n}{2(\mathtt{k}^+ -\delta)}(\frac{1}{m}-\frac{1}{2})+ \frac{\mathtt{k}^-}{2(\mathtt{k}^+ -\delta)}} \big(\|u_0\|_{L^m \cap H^{\mathtt{k}^+}}+ \|u_1\|_{L^m \cap L^2}\big), \\
\big\||D|^{\mathtt{k}^+} u(t,\cdot)\big\|_{L^2}& \lesssim (1+t)^{-\frac{n}{2(\mathtt{k}^+ -\delta)}(\frac{1}{m}-\frac{1}{2})- \frac{\mathtt{k}^+- \mathtt{k}^-}{2(\mathtt{k}^+ -\delta)}} \big(\|u_0\|_{L^m \cap H^{\mathtt{k}^+}}+ \|u_1\|_{L^m \cap L^2}\big), \\
\|u_t(t,\cdot)\|_{L^2}& \lesssim (1+t)^{-\frac{n}{2(\mathtt{k}^+ -\delta)}(\frac{1}{m}-\frac{1}{2})- \frac{\sigma- \mathtt{k}^-}{\mathtt{k}^+ -\delta}} \big(\|u_0\|_{L^m \cap H^{\mathtt{k}^+}}+ \|u_1\|_{L^m \cap L^2}\big).
\end{align*}
\edl
We are going to prove the following main result.
\begin{dl}[\textbf{Blow-up}] \label{dloptimal9.1}
Let $\sigma \ge 1$, $\delta \in [0,\sigma)$ and $n > 2\mathtt{k}^-$. We assume that we choose the initial data $u_0=0$ and $u_1 \in L^1$ satisfying the following relation:
\begin{equation} \label{optimal9.1}
\int_{\R^n} u_1(x)dx >0.
\end{equation}
Moreover, we suppose the condition
\begin{equation} \label{optimal9.2}
p \in \Big(1, 1+ \frac{2\sigma}{n- \mathtt{k}^-}\Big).
\end{equation}
Then, there is no global (in time) Sobolev solution $u \in C\big([0,\infty),H^{2\sigma}\big)$ to (\ref{pt9.1}).
\end{dl}

\bnx \label{remark9.1.2}
\fontshape{n}
\selectfont
We want to underline that the lifespan $T_\e$ of Sobolev solutions to given data $(0,\e u_1)$ for any small positive constant $\e$ in Theorem \ref{dloptimal9.1} can be estimated as follows:
\begin{equation} \label{lifetime9.1}
T_\e \le C\e^{-\frac{(2\sigma- \mathtt{k}^-)(p-1)}{2\sigma- (n- \mathtt{k}^-)(p-1)}} \quad \text{ with } C>0.
\end{equation}
\enx

\begin{nx}
\fontshape{n}
\selectfont
If we choose $m=1$ in Theorem \ref{dlexistence9.1}, then from Theorem \ref{dloptimal9.1} it is clear that the critical exponent $p_{crit}$ is given by
$$ p_{crit}(n)=
1+\frac{2\sigma}{n-2\delta} \quad \text{ if }\delta \in \Big[0,\frac{\sigma}{2}\Big] \text{ and }4\delta <n \leq 4\sigma.$$
However, in the case $\delta \in (\frac{\sigma}{2},\sigma)$ there appears a gap between the exponents given by $1+\frac{2\delta+\sigma}{n-\sigma}$ from Theorem \ref{dlexistence9.1} and  $1+\frac{2\sigma}{n-\sigma}$ from Theorem  \ref{dloptimal9.1} for $2\sigma <n \leq 8\delta$.
\end{nx}

\section{Preliminaries} \label{Sec.Pre}
In this section, we collect some preliminary knowledge needed in our proofs.
\bdn[\cite{Kwanicki,Silvestre}]
\fontshape{n}
\selectfont
Let $s \in (0,1)$. Let $X$ be a suitable set of functions defined on $\R^n$. Then, the fractional Laplacian $(-\Delta)^s$ in $\R^n$ is a non-local operator given by
$$ (-\Delta)^s: \,\,v \in X  \to (-\Delta)^s v(x):= C_{n,s}\,\, p.v.\int_{\R^n}\frac{v(x)- v(y)}{|x-y|^{n+2s}}dy $$
as long as the right-hand side exists, where $p.v.$ stands for Cauchy's principal value, $C_{n,s}:= \frac{4^s \Gamma(\frac{n}{2}+s)}{\pi^{\frac{n}{2}}\Gamma(-s)}$ is a normalization constant and $\Gamma$ denotes the Gamma function.
\label{def1}
\edn

\bbd \label{lemma2.1}
Let $\big< x\big>=(1+|x|^2)^{\frac{1}{2}}$ and $q>0$. Then, the following estimate holds for any multi-index $\alpha$ satisfying $|\alpha|\ge 1$:
$$ \big|\partial_x^\alpha \big< x\big>^{-q}\big| \lesssim \big< x\big>^{-q-|\alpha|}. $$
\ebd
\begin{proof}
First, we recall the following formula of derivatives of composed functions for $|\alpha|\ge 1$:
$$ \partial_x^\alpha h\big(f(x)\big)= \sum_{k=1}^{|\alpha|}h^{(k)} \big(f(x)\big)\left(\sum_{\substack{\gamma_1+\cdots+\gamma_k \le \alpha\\ |\gamma_1|+\cdots+|\gamma_k|= |\alpha|,\, |\gamma_i|\ge 1}}\big(\partial_x^{\gamma_1} f(x)\big) \cdots \big(\partial_x^{\gamma_k} f(x)\big)\right), $$
where $h=h(z)$ and $h^{(k)}(z)=\frac{d^k h(z)}{dz^k}$. Applying this formula with $h(z)= z^{-\frac{q}{2}}$ and $f(x)= 1+|x|^2$ we obtain
\begin{align*}
\big|\partial_x^\alpha \big< x\big>^{-q}\big|&\le \sum_{k=1}^{|\alpha|} (1+|x|^2)^{-\frac{q}{2}-k}\left(\sum_{\substack{\gamma_1+\cdots+\gamma_k \le \alpha\\ |\gamma_1|+\cdots+|\gamma_k|= |\alpha|,\, |\gamma_i|\ge 1}}\big|\partial_x^{\gamma_1} (1+|x|^2)\big| \cdots \big|\partial_x^{\gamma_k} (1+|x|^2)\big|\right) \\
&\le C_1\sum_{k=1}^{|\alpha|} (1+|x|^2)^{-\frac{q}{2}-k}
\begin{cases}
1 &\quad \text{ if } 0 \le |x| \le 1, \\
\left(\displaystyle\sum_{\substack{\gamma_1+\cdots+\gamma_k \le \alpha\\ |\gamma_1|+\cdots+|\gamma_k|= |\alpha|,\, |\gamma_i|\ge 1}}|x|^{2-|\gamma_1|} \cdots |x|^{2-|\gamma_k|}\right) &\quad \text{ if } |x| \ge 1,
\end{cases} \\
&\le C_2\sum_{k=1}^{|\alpha|} (1+|x|^2)^{-\frac{q}{2}-k}
\begin{cases}
1 &\quad \text{ if }\quad 0 \le |x| \le 1, \\
|x|^{2k-|\alpha|} &\quad \text{ if }\quad |x| \ge 1,
\end{cases} \\
&\le \begin{cases}
C_2|\alpha|\big<x\big>^{-q-2}&\quad \text{ if }\quad 0 \le |x| \le 1, \\
C_2|\alpha|\big< x\big>^{-q}|x|^{-|\alpha|}&\quad \text{ if }\quad |x| \ge 1,
\end{cases}
\end{align*}
where $C_1$ and $C_2$ are some suitable constants. This completes the proof.
\end{proof}

\bbd \label{bd9.3.1}
Let $m \in \Z$, $s \in (0,1)$ and $\gamma:= m+ s$. If $v \in H^{2\gamma}(\R^n)$, then it holds
\begin{equation*}
(-\Delta)^{\gamma}v(x)= (-\Delta)^{m}\big((-\Delta)^{s}v(x)\big)= (-\Delta)^{s}\big((-\Delta)^{m}v(x)\big).
\end{equation*}
\ebd
\noindent One can find the proof of Lemma \ref{bd9.3.1} in Remark $3.2$ in \cite{Abatangelo}.

\bbd \label{lemma2.2}
Let $m \in \Z$, $s \in (0,1)$ and $\gamma:= m+ s$. Let $\big< x\big>=(1+|x|^2)^{\frac{1}{2}}$ and $q>0$. Then, the following estimates hold for all $x \in \R^n$:
\begin{equation}
\big|(-\Delta)^\gamma \big< x\big>^{-q}\big| \lesssim
\begin{cases}
\big< x\big>^{-q-2\gamma} &\quad \text{ if }\quad 0< q+2m< n, \\
\big< x\big>^{-n-2s}\log(e+|x|) &\quad \text{ if }\quad q+2m= n, \\
\big< x\big>^{-n-2s} &\quad \text{ if }\quad q+2m> n.
\end{cases} \label{lemma2.2.Estimate}
\end{equation}
\ebd
\begin{proof}
We follow ideas from the proof of Lemma $1.5$ in \cite{Fujiwara} devoting to the case $m=0$ and $s=\frac{1}{2}$, that is, the case $\gamma= \frac{1}{2}$ is generalized to any fractional number $\gamma>0$. To do this, for any $s \in (0,1)$ we shall divide the proof into two cases: $m=0$ and $m\ge 1$. \smallskip

\noindent Let us consider the first case $m=0$. Denoting by $\psi= \psi(x):= \big< x\big>^{-q}$ we write $(-\Delta)^s \big< x\big>^{-q}= (-\Delta)^s (\psi)(x)$. According to Definition \ref{def1} of fractional Laplacian as a singular integral operator, we have
$$(-\Delta)^s (\psi)(x):= C_{n,\delta}\,\, p.v.\int_{\R^n}\frac{\psi(x)- \psi(y)}{|x-y|^{n+ 2s}}dy. $$
A standard change of variables leads to
\begin{align*}
(-\Delta)^s (\psi)(x)&= -\frac{C_{n,s}}{2}\,\, p.v.\int_{\R^n}\frac{\psi(x+y)+ \psi(x-y)- 2\psi(x)}{|y|^{n+ 2s}}dy \\
&= -\frac{C_{n,s}}{2}\lim_{\e \to 0+}\int_{\e\le |y|\le 1}\frac{\psi(x+y)+ \psi(x-y)- 2\psi(x)}{|y|^{n+ 2s}}dy \\
&\quad- \frac{C_{n,s}}{2}\int_{|y|\ge 1}\frac{\psi(x+y)+ \psi(x-y)- 2\psi(x)}{|y|^{n+ 2s}}dy.
\end{align*}
To deal with the first integral, after using a second order Taylor expansion for $\psi$ we arrive at
$$ \frac{|\psi(x+y)+ \psi(x-y)- 2\psi(x)|}{|y|^{n+2s}} \lesssim \frac{\|\partial_x^2 \psi\|_{L^\ity}}{|y|^{n+ 2s- 2}}. $$
Thanks to the above estimate and $s\in (0,1)$, we may remove the principal value of the integral at the origin to conclude
$$ (-\Delta)^s (\psi)(x)= -\frac{C_{n,s}}{2} \int_{\R^n}\frac{\psi(x+y)+ \psi(x-y)- 2\psi(x)}{|y|^{n+ 2s}}dy. $$
To prove the desired estimates, we shall divide our considerations into two cases. In the first subcase $\{x:|x|\le 1\}$, we can proceed as follows:
\begin{align*}
\big|(-\Delta)^s (\psi)(x)\big|&\lesssim \int_{|y|\le 1}\frac{|\psi(x+y)+ \psi(x-y)- 2\psi(x)|}{|y|^{n+ 2s}}dy+ \int_{|y|\ge 1}\frac{|\psi(x+y)+ \psi(x-y)- 2\psi(x)|}{|y|^{n+ 2s}}dy \\
&\lesssim \|\partial_x^2 \psi\|_{L^\ity} \int_{|y|\le 1}\frac{1}{|y|^{n+ 2s- 2}}dy+ \|\psi\|_{L^\ity} \int_{|y|\ge 1}\frac{1}{|y|^{n+ 2s}}dy.
\end{align*}
Due to the boundedness of the above two integrals, it follows immediately
\begin{equation}
\big|(-\Delta)^s (\psi)(x)\big| \lesssim 1 \quad \text{ for } |x|\le 1. \label{lem2.2.1}
\end{equation}
In order to deal with the second subcase $\{x:|x|\ge 1\}$, we can re-write
\begin{align}
(-\Delta)^s (\psi)(x)&= -\frac{C_{n,s}}{2}\int_{|y|\ge 2|x|} \frac{\psi(x+y)+ \psi(x-y)- 2\psi(x)}{|y|^{n+ 2s}}dy \nonumber \\
&\quad- \frac{C_{n,s}}{2}\int_{\frac{1}{2}|x|\le |y|\le 2|x|} \frac{\psi(x+y)+ \psi(x-y)- 2\psi(x)}{|y|^{n+ 2s}}dy \nonumber \\
&\quad- \frac{C_{n,s}}{2}\int_{|y|\le \frac{1}{2}|x|} \frac{\psi(x+y)+ \psi(x-y)- 2\psi(x)}{|y|^{n+ 2s}}dy. \label{lem2.2.2}
\end{align}
For the first integral, we notice that the relations $|x+y|\ge |y|-|x| \ge |x|$ and $|x-y|\ge |y|-|x| \ge |x|$ hold for $|y|\ge 2|x|$. Since $\psi$ is a decreasing function, we obtain the following estimate:
\begin{align}
&\Big|\int_{|y|\ge 2|x|} \frac{\psi(x+y)+ \psi(x-y)- 2\psi(x)}{|y|^{n+ 2s}}dy\Big| \nonumber \\
&\qquad \le 4|\psi(x)| \int_{|y|\ge 2|x|} \frac{1}{|y|^{n+ 2s}}dy \lesssim \big< x\big>^{-q} \int_{|y|\ge 2|x|} \frac{1}{|y|^{1+ 2s}}d|y| \nonumber \\
&\qquad \lesssim \big< x\big>^{-q} |x|^{-2s}\lesssim \big< x\big>^{-q-2s} \qquad \big(\text{due to } |x| \approx \big< x\big> \text{ for } |x|\ge 1\big). \label{lem2.2.3}
\end{align}
It is clear that $|y| \approx |x|$ in the second integral domain. Moreover, it follows
\begin{align}
\Big\{y: \frac{1}{2}|x|\le |y|\le 2|x|\Big\} &\subset \big\{y: |x+y|\le 3|x|\big\}, \label{lem2.2.31}\\
\Big\{y: \frac{1}{2}|x|\le |y|\le 2|x|\Big\} &\subset \big\{y: |x-y|\le 3|x|\big\}. \label{lem2.2.32}
\end{align}
For this reason, we arrive at
\begin{align}
&\Big|\int_{\frac{1}{2}|x|\le |y|\le 2|x|} \frac{\psi(x+y)+ \psi(x-y)- 2\psi(x)}{|y|^{n+ 2s}}dy\Big| \nonumber \\
&\qquad \lesssim |x|^{-n-2s}\Big(\int_{|x+y|\le 3|x|} \psi(x+y)dy+ \int_{|x-y|\le 3|x|} \psi(x-y)dy+ \psi(x)\int_{\frac{1}{2}|x|\le |y|\le 2|x|} 1 dy\Big) \nonumber \\
&\qquad \lesssim |x|^{-n-2s}\Big(\int_{|x+y|\le 3|x|} \psi(x+y)dy+ \big< x\big>^{-q}|x|^n\Big), \label{lem2.2.4}
\end{align}
where we used the relation \[ \int_{|x+y|\le 3|x|} \psi(x+y)dy= \int_{|x-y|\le 3|x|} \psi(x-y)dy.\] By the change of variables $r=|x+y|$, we apply the inequality $1+r^2 \ge \frac{(1+r)^2}{2}$ to get
\begin{align}
\int_{|x+y|\le 3|x|} \psi(x+y)dy&\lesssim \int_{r \le 3|x|} (1+r^2)^{-\frac{q}{2}}\,r^{n-1}dr \lesssim \int_{r \le 3|x|} (1+r)^{n-q-1}dr \nonumber \\
&\lesssim \begin{cases}
(1+3|x|)^{n-q} &\quad \text{ if }\quad 0< q< n, \\
\log(e+3|x|) &\quad \text{ if }\quad q= n, \\
1 &\quad \text{ if }\quad q> n.
\end{cases} \label{lem2.2.5}
\end{align}
By $|x| \approx \big< x\big>$ for $|x|\ge 1$, combining (\ref{lem2.2.4}) and (\ref{lem2.2.5}) leads to
\begin{equation}
\Big|\int_{\frac{1}{2}|x|\le |y|\le 2|x|} \frac{\psi(x+y)+ \psi(x-y)- 2\psi(x)}{|y|^{n+ 2s}}dy\Big| \lesssim
\begin{cases}
\big< x\big>^{-q-2s} & \text{ if }\quad 0< q< n, \\
\big< x\big>^{-n-2s}\log(e+3|x|) & \text{ if }\quad q= n, \\
\big< x\big>^{-n-2s}& \text{ if }\quad q> n.
\end{cases}
\label{lem2.2.6}
\end{equation}
For the third integral in (\ref{lem2.2.2}), using again the second order Taylor expansion for $\psi$ we obtain
\begin{align}
&\Big|\int_{|y|\le \frac{1}{2}|x|} \frac{\psi(x+y)+ \psi(x-y)- 2\psi(x)}{|y|^{n+ 2s}}dy\Big| \nonumber \\
&\qquad \le \int_{|y|\le \frac{1}{2}|x|} \frac{|\psi(x+y)+ \psi(x-y)- 2\psi(x)|}{|y|^{n+ 2s}}dy \lesssim \int_{|y|\le \frac{1}{2}|x|} \max_{\theta \in [0,1]}\big|\partial_x^2 \psi(x\pm \theta y)\big|\frac{1}{|y|^{n+ 2s- 2}}dy \nonumber \\
&\qquad \lesssim \int_{|y|\le \frac{1}{2}|x|} \max_{\theta \in [0,1]}\big< x\pm \theta y\big>^{-q-2}\frac{1}{|y|^{n+ 2s- 2}}dy \lesssim \big< x\big>^{-q-2}\int_{|y|\le \frac{1}{2}|x|} |y|^{1-2s}d|y| \lesssim \big< x\big>^{-q-2s}. \label{lem2.2.7}
\end{align}
Here we used the relation $|x\pm \theta y| \ge |x|- \theta |y|\ge |x|- \frac{1}{2}|x|= \frac{1}{2}|x|$. From (\ref{lem2.2.2}), (\ref{lem2.2.3}), (\ref{lem2.2.6}) and (\ref{lem2.2.7}) we arrive at the following estimates for $|x|\ge 1$:
\begin{equation}
\big|(-\Delta)^s (\psi)(x)\big| \lesssim
\begin{cases}
\big< x\big>^{-q-2s} &\quad \text{ if }\quad 0< q< n, \\
\big< x\big>^{-n-2s}\log(e+3|x|) &\quad \text{ if }\quad q= n, \\
\big< x\big>^{-n-2s}&\quad \text{ if }\quad q> n.
\end{cases}
\label{lem2.2.8}
\end{equation}
Finally, combining (\ref{lem2.2.1}) and (\ref{lem2.2.8}) we may conclude all desired estimates for $m=0$. \smallskip

\noindent Next let us turn to the second case $m\ge 1$. First, a straight-forward calculation gives the following relation:
\begin{equation}
-\Delta \big< x\big>^{-r}= r\Big((n-r-2)\big< x\big>^{-r-2} + (r+2)\big< x\big>^{-r-4}\Big) \quad \text{ for any } r>0. \label{RepresentationFormula^1}
\end{equation}
By induction argument, carrying out $m$ steps of (\ref{RepresentationFormula^1}) we obtain the following formula for any $m\ge 1$:
\begin{align}
(-\Delta)^m \big< x\big>^{-q}&= (-1)^m \prod_{j=0}^{m-1}(q+2j)\Big(\prod_{j=1}^{m}(-n+q+2j)\big< x\big>^{-q-2m} \nonumber \\
&\hspace{4cm}- C^1_m \prod_{j=2}^{m}(-n+q+2j)(q+2m)\big< x\big>^{-q-2m-2} \nonumber \\
&\hspace{4cm}+ C^2_m \prod_{j=3}^{m}(-n+q+2j)(q+2m)(q+2m+2)\big< x\big>^{-q-2m-4} \nonumber \\
&\hspace{4cm}+\cdots+ (-1)^m \prod_{j=0}^{m-1}(q+2m+2j)\big< x\big>^{-q-4m}\Big). \label{RepresentationFormula^m}
\end{align}
Then, thanks to Lemma \ref{bd9.3.1}, we derive
\begin{align}
(-\Delta)^\gamma \big< x\big>^{-q}&= (-\Delta)^s \big((-\Delta)^m \big< x\big>^{-q}\big) \nonumber \\
&= (-1)^m \prod_{j=0}^{m-1}(q+2j)\Big(\prod_{j=1}^{m}(-n+q+2j)\, (-\Delta)^s \big< x\big>^{-q-2m} \nonumber \\
&\hspace{3cm}- C^1_m \prod_{j=2}^{m}(-n+q+2j)(q+2m)\, (-\Delta)^s \big< x\big>^{-q-2m-2}  \nonumber \\
&\hspace{3cm}+ C^2_m \prod_{j=3}^{m}(-n+q+2j)(q+2m)(q+2m+2)\, (-\Delta)^s \big< x\big>^{-q-2m-4}  \nonumber \\
&\hspace{3cm}+\cdots+ (-1)^m \prod_{j=0}^{m-1}(q+2m+2j)\, (-\Delta)^s \big< x\big>^{-q-4m}\Big). \label{lem2.2.9}
\end{align}
For this reason, in order to conclude the desired estimates, we only indicate the following estimates for $k=0,\cdots,m$:
\begin{equation}
\big|(-\Delta)^s \big< x\big>^{-q-2(m+k)}\big| \lesssim
\begin{cases}
\big< x\big>^{-q-2\gamma} &\quad \text{ if }\quad 0< q+2m< n, \\
\big< x\big>^{-n-2s}\log(e+|x|) &\quad \text{ if }\quad q+2m= n, \\
\big< x\big>^{-n-2s} &\quad \text{ if }\quad q+2m> n.
\end{cases}
\label{lem2.2.10}
\end{equation}
Indeed, substituting $q$ by $q+2(m+k)$ with $k=0,\cdots,m$ and $\gamma= s$ into (\ref{lemma2.2.Estimate}) leads to
$$ \big|(-\Delta)^s \big< x\big>^{-q-2(m+k)}\big| \lesssim
\begin{cases}
\big< x\big>^{-q-2\gamma} &\quad \text{ if }\quad 0< q+2(m+k)< n, \\
\big< x\big>^{-n-2s}\log(e+|x|) &\quad \text{ if }\quad q+2(m+k)= n, \\
\big< x\big>^{-n-2s} &\quad \text{ if }\quad q+2(m+k)> n.
\end{cases}$$
From these estimates, it follows immediately (\ref{lem2.2.10}) to conclude (\ref{lemma2.2.Estimate}) for any $m\ge 1$. Summarizing, the proof of Lemma \ref{lemma2.2} is completed.
\end{proof}

\bbd \label{lemma2.3}
Let $s\in (0,1)$. Let $\psi$ be a smooth function satisfying $\partial_x^2 \psi\in L^\ity$. For any $R>0$, let $\psi_R$ be a function defined by
$$ \psi_R(x):= \psi\big(R^{-1} x\big) $$
for all $x \in \R^n$. Then, $(-\Delta)^s (\psi_R)$ satisfies the following scaling properties for all $x \in \R^n$:
$$(-\Delta)^s (\psi_R)(x)= R^{-2s}\big((-\Delta)^s \psi \big)\big(R^{-1} x\big). $$
\ebd
\begin{proof}
Thanks to the assumption $\partial_x^2 \psi\in L^\ity$, following the proof of Lemma \ref{lemma2.2} we may remove the principal value of the integral at the origin to conclude
\begin{align*}
(-\Delta)^s (\psi_R)(x) &= -\frac{C_{n,s}}{2} \int_{\R^n}\frac{\psi_R(x+y)+ \psi_R(x-y)- 2\psi_R(x)}{|y|^{n+ 2s}}dy \\
&= -\frac{C_{n,s}}{2} R^{-2s} \int_{\R^n}\frac{\psi\big(R^{-1} x+ R^{-1} y\big)+ \psi\big(R^{-1} x- R^{-1} y\big)- 2\psi\big(R^{-1} x\big)}{|R^{-1} y|^{n+ 2s}}d(R^{-1} y) \\
&= R^{-2s}\big((-\Delta)^s \psi\big)\big(R^{-1} x\big).
\end{align*}
This completes the proof.
\end{proof}

\section{Proof of the blow-up result} \label{Sec.Proof}
\noindent We divide the proof of Theorem \ref{dloptimal9.1} into several cases.

\subsection{The case that both parameters $\sigma$ and $\delta$ are integers} \label{Sec9.1}
The proof of this case can be found in the paper \cite{DabbiccoEbert}.

\subsection{The case that the parameter $\sigma$ is integer and the parameter $\delta$ is fractional from $(0,1)$} \label{Sec9.2}
\begin{proof}
First, we introduce the function $\varphi=\varphi(|x|):=\big< x\big>^{-n-2\delta}$ and the function $\eta= \eta(t)$ having the following properties:
\begin{align}
&1.\quad \eta \in \mathcal{C}_0^\ity([0,\ity)) \text{ and }
\eta(t)=\begin{cases}
1 &\quad \text{ for }0 \le t \le \frac{1}{2}, \\
\text{decreasing } &\quad \text{ for }\frac{1}{2} \le t \le 1, \\
0 &\quad \text{ for }t \ge 1,
\end{cases} \nonumber \\
&2.\quad \eta^{-\frac{p'}{p}}(t)\big(|\eta'(t)|^{p'}+|\eta''(t)|^{p'}\big) \le C \quad \text{ for any } t \in \Big[\frac{1}{2},1\Big], \label{condition9.2.1}
\end{align}
where $p'$ is the conjugate of $p>1$. Let $R$ be a large parameter in $[0,\ity)$. We define the following test function:
$$ \phi_R(t,x):= \eta_R(t) \varphi_R(x), $$
where $\eta_R(t):= \eta\big(R^{-\alpha}t\big)$ and $\varphi_R(x):= \varphi\big(R^{-1}x\big)$ with a fixed parameter $\alpha:= 2\sigma- \mathtt{k}^-$. We define the functionals
$$ I_R:= \int_0^{\ity}\int_{\R^n}|u(t,x)|^p \phi_R(t,x)\,dxdt= \int_0^{R^{\alpha}}\int_{\R^n}|u(t,x)|^p \phi_R(t,x)\,dxdt $$
and
$$ I_{R,t}:= \int_{\frac{R^\alpha}{2}}^{R^{\alpha}}\int_{\R^n}|u(t,x)|^p \phi_R(t,x)\,dxdt. $$
Let us assume that $u= u(t,x)$ is a global (in time) Sobolev solution from $C\big([0,\infty),H^{2\sigma}\big)$ to (\ref{pt9.1}). After multiplying the equation (\ref{pt9.1}) by $\phi_R=\phi_R(t,x)$, we carry out partial integration to derive
\begin{align}
0\le I_R &= -\int_{\R^n} u_1(x)\varphi_R(x)\,dx + \int_{\frac{R^\alpha}{2}}^{R^{\alpha}}\int_{\R^n}u(t,x) \partial_t^2 \eta_R(t) \varphi_R(x)\,dxdt \nonumber \\
&\quad + \int_0^{\ity}\int_{\R^n} \eta_R(t) \varphi_R(x)\, (-\Delta)^{\sigma} u(t,x)\,dxdt- \int_{\frac{R^\alpha}{2}}^{R^\alpha}\int_{\R^n} \partial_t \eta_R(t) \varphi_R(x)\,(-\Delta)^{\delta} u(t,x)\,dxdt \nonumber \\
&=: -\int_{\R^n} u_1(x)\varphi_R(x)\,dx+ J_1+ J_2- J_3. \label{t9.2.1}
\end{align}
Applying H\"{o}lder's inequality with $\frac{1}{p}+\frac{1}{p'}=1$ we may estimate as follows:
\begin{align*}
|J_1| &\le \int_{\frac{R^\alpha}{2}}^{R^{\alpha}}\int_{\R^n} |u(t,x)|\, \big|\partial_t^2 \eta_R(t)\big| \varphi_R(x) \, dxdt \\
&\lesssim \Big(\int_{\frac{R^\alpha}{2}}^{R^{\alpha}}\int_{\R^n} \Big|u(t,x)\phi^{\frac{1}{p}}_R(t,x)\Big|^p \,dxdt\Big)^{\frac{1}{p}} \Big(\int_{\frac{R^\alpha}{2}}^{R^{\alpha}}\int_{\R^n} \Big|\phi^{-\frac{1}{p}}_R(t,x) \partial_t^2 \eta_R(t) \varphi_R(x)\Big|^{p'}\, dxdt\Big)^{\frac{1}{p'}} \\
&\lesssim I_{R,t}^{\frac{1}{p}}\, \Big(\int_{\frac{R^\alpha}{2}}^{R^{\alpha}}\int_{\R^n} \eta_R^{-\frac{p'}{p}}(t) \big|\partial_t^2 \eta_R(t)\big|^{p'} \varphi_R(x)\, dxdt\Big)^{\frac{1}{p'}}.
\end{align*}
By the change of variables $\tilde{t}:= R^{-\alpha}t$ and $\tilde{x}:= R^{-1}x$, a straight-forward calculation gives
\begin{equation}
|J_1| \lesssim I_{R,t}^{\frac{1}{p}}\, R^{-2\alpha+ \frac{n+\alpha}{p'}}\Big(\int_{\R^n} \big< \tilde{x}\big>^{-n-2\delta}\, d\tilde{x}\Big)^{\frac{1}{p'}}.
\label{t9.2.2}
\end{equation}
Here we used $\partial_t^2 \eta_R(t)= R^{-2\alpha}\eta''(\tilde{t})$ and the assumption (\ref{condition9.2.1}). Now let us turn to estimate $J_2$ and $J_3$. First, we notice that by Parseval-Plancherel formula it holds:
$$ \int_{\R^n}v_1(x)\,(-\Delta)^{\gamma}v_2(x) dx= \int_{\R^n}|\xi|^{2\gamma}\widehat{v}_1(\xi)\,\widehat{v}_2(\xi) d\xi= \int_{\R^n}v_2(x)\,(-\Delta)^{\gamma}v_1(x) dx, $$
for any $\gamma > 0$ and $v_1,\,v_2 \in H^{2\gamma}$. Here $\widehat{v}_j=\widehat{v}_j(\xi)$ stands for the Fourier transform with respect to the spatial variables of $v_j=v_j(x)$,  $j=1,2$. Using $\varphi_R \in H^{2\sigma}$ and $u \in C\big([0,\infty),H^{2\sigma}\big)$ we may conclude 
\begin{align*}
\int_{\R^n} \varphi_R(x)\, (-\Delta)^{\sigma} u(t,x)\,dx&= \int_{\R^n} u(t,x)\, (-\Delta)^{\sigma}\varphi_R(x) \,dx, \\
\int_{\R^n} \varphi_R(x)\,(-\Delta)^{\delta} u(t,x)\,dx&= \int_{\R^n} u(t,x)\,(-\Delta)^{\delta}\varphi_R(x)\, dx,
\end{align*}
that is,
\begin{align*}
\big(\varphi_R,(-\Delta)^{\sigma} u(t,\cdot)\big)_{L^2}&= \big((-\Delta)^{\sigma}\varphi_R,\,u(t,\cdot)\big)_{L^2}, \\
\big(\varphi_R,(-\Delta)^{\delta} u(t,\cdot)\big)_{L^2}&= \big((-\Delta)^{\delta}\varphi_R,\,u(t,\cdot)\big)_{L^2}.
\end{align*}
We can see that under the assumptions both scalar products are well defined. Hence, we obtain
$$J_2= \int_0^{\ity}\int_{\R^n} \eta_R(t) \varphi_R(x)\, (-\Delta)^{\sigma} u(t,x)\,dxdt= \int_0^{\ity}\int_{\R^n} \eta_R(t) u(t,x)\, (-\Delta)^{\sigma}\varphi_R(x) \,dxdt, $$
and
$$J_3= \int_{\frac{R^\alpha}{2}}^{R^\alpha}\int_{\R^n} \partial_t \eta_R(t) \varphi_R(x)\,(-\Delta)^{\delta} u(t,x)\,dxdt= \int_{\frac{R^\alpha}{2}}^{R^\alpha}\int_{\R^n} \partial_t \eta_R(t) u(t,x)\,(-\Delta)^{\delta}\varphi_R(x)\, dxdt. $$
Applying H\"{o}lder's inequality again as we estimated $J_1$ leads to
$$|J_2|\le I_R^{\frac{1}{p}}\, \Big(\int_0^{R^{\alpha}}\int_{\R^n} \eta_R(t) \varphi^{-\frac{p'}{p}}_R(x)\, \big|(-\Delta)^{\sigma}\varphi_R(x)\big|^{p'} \, dxdt\Big)^{\frac{1}{p'}}, $$
and
$$|J_3|\le I_{R,t}^{\frac{1}{p}}\, \Big(\int_{\frac{R^\alpha}{2}}^{R^{\alpha}}\int_{\R^n} \eta^{-\frac{p'}{p}}_R(t) \big|\partial_t \eta_R(t)\big|^{p'} \varphi^{\frac{-p'}{p}}_R(x)\, \big|(-\Delta)^{\delta}\varphi_R(x)\big|^{p'} \, dxdt\Big)^{\frac{1}{p'}}. $$
In order to control the above two integrals, the key tools rely on results from Lemmas \ref{lemma2.1}, \ref{lemma2.2} and \ref{lemma2.3}. Namely, at first carrying out the change of variables $\tilde{t}:= R^{-\alpha}t$ and $\tilde{x}:= R^{-1}x$ we arrive at
\begin{align*}
|J_2| &\lesssim I_R^{\frac{1}{p}}\, R^{-2\sigma+ \frac{n+\alpha}{p'}}\Big(\int_0^{1}\int_{\R^n} \eta(\tilde{t}) \varphi^{-\frac{p'}{p}}(\tilde{x})\, \big|(-\Delta)^{\sigma}(\varphi)(\tilde{x})\big|^{p'}\, d\tilde{x}d\tilde{t}\Big)^{\frac{1}{p'}} \\
&\lesssim I_R^{\frac{1}{p}}\, R^{-2\sigma+ \frac{n+\alpha}{p'}}\Big(\int_{\R^n} \varphi^{-\frac{p'}{p}}(\tilde{x})\, \big|(-\Delta)^{\sigma}(\varphi)(\tilde{x})\big|^{p'}\, d\tilde{x}\Big)^{\frac{1}{p'}},
\end{align*}
where we note ($\sigma$ is an integer) that $(-\Delta)^{\sigma}\varphi_R(x)= R^{-2\sigma}(-\Delta)^{\sigma}\varphi(\tilde{x}).$ Using Lemma \ref{lemma2.1} implies the following estimate:
\begin{equation}
|J_2|\lesssim I_R^{\frac{1}{p}}\, R^{-2\sigma+ \frac{n+\alpha}{p'}}\Big(\int_{\R^n} \big< \tilde{x}\big>^{-n-2\delta-2\sigma p'}\, d\tilde{x}\Big)^{\frac{1}{p'}}.
 \label{t9.2.3}
\end{equation}
Next carrying out again the change of variables $\tilde{t}:= R^{-\alpha}t$ and $\tilde{x}:= R^{-1}x$ and employing Lemma \ref{lemma2.3} we can proceed $J_3$ as follows:
\begin{align*}
|J_3| &\lesssim I_{R,t}^{\frac{1}{p}}\, R^{-2\delta-\alpha+ \frac{n+\alpha}{p'}}\Big(\int_{\frac{1}{2}}^{1}\int_{\R^n} \eta^{-\frac{p'}{p}}(\tilde{t}) \big|\eta'(\tilde{t})\big|^{p'} \varphi^{-\frac{p'}{p}}(\tilde{x})\, \big|(-\Delta)^{\delta}(\varphi)(\tilde{x})\big|^{p'}\, d\tilde{x}d\tilde{t}\Big)^{\frac{1}{p'}} \\
&\lesssim I_{R,t}^{\frac{1}{p}}\, R^{-2\delta-\alpha+ \frac{n+\alpha}{p'}}\Big(\int_{\R^n} \varphi^{-\frac{p'}{p}}(\tilde{x})\, \big|(-\Delta)^{\delta}(\varphi)(\tilde{x})\big|^{p'}\, d\tilde{x}\Big)^{\frac{1}{p'}}.
\end{align*}
Here we used $\partial_t \eta_R(t)= R^{-\alpha}\eta'(\tilde{t})$ and the assumption (\ref{condition9.2.1}). To deal with the last integral, we apply Lemma \ref{lemma2.2} with $q=n+2\delta$ and $\gamma=\delta$, that is, $m=0$ and $s=\delta$ to get
\begin{equation}
|J_3| \lesssim I_{R,t}^{\frac{1}{p}} R^{-2\delta-\alpha+ \frac{n+\alpha}{p'}}\Big(\int_{\R^n} \big< \tilde{x}\big>^{-n-2\delta}\, d\tilde{x}\Big)^{\frac{1}{p'}}.
\label{t9.2.4}
\end{equation}
Because of assumption (\ref{optimal9.1}), there exists a sufficiently large constant $R_0> 0$ such that it holds
\begin{equation}
\int_{\R^n} u_1(x) \varphi_R(x)\, dx >0
 \label{t9.2.5}
\end{equation}
for all $R > R_0$. Combining the estimates from (\ref{t9.2.1}) to (\ref{t9.2.5}) we may arrive at
\begin{align}
0< \int_{\R^n} u_1(x) \varphi_R(x)\, dx &\lesssim I_{R,t}^{\frac{1}{p}} \big(R^{-2\alpha+ \frac{n+\alpha}{p'}}+ R^{-\alpha- 2\delta+ \frac{n+\alpha}{p'}}\big)+ I_R^{\frac{1}{p}}\, R^{-2\sigma+ \frac{n+\alpha}{p'}}- I_R \label{t9.2.6} \\
&\lesssim I_R^{\frac{1}{p}}R^{-2\sigma+ \frac{n+\alpha}{p'}}- I_R \label{t9.2.7}
\end{align}
for all $R > R_0$. Moreover, applying the inequality
$$ A\,y^\gamma- y \le A^{\frac{1}{1-\gamma}} \quad \text{ for any } A>0,\, y \ge 0 \text{ and } 0< \gamma< 1 $$
leads to
\begin{equation}
0< \int_{\R^n} u_1(x)\varphi_R(x)dx \lesssim R^{-2\sigma p'+ n+\alpha}
\label{t9.2.8}
\end{equation}
for all $R > R_0$. It is clear that the assumption (\ref{optimal9.2}) is equivalent to $-2\sigma p'+ n+\alpha< 0$. For this reason, letting $R \to \ity$ in (\ref{t9.2.8}) we obtain
$$ \int_{\R^n} u_1(x)\,dx= 0. $$
This is a contradiction to the assumption (\ref{optimal9.1}). Summarizing, the proof is completed.
\end{proof}

Let us now consider the case of subcritical exponent to explain the estimate for lifespan $T_\e$ of solutions in Remark \ref{remark9.1.2}. We assume that $u= u(t,x)$ is a local (in time) Sobolev solution to (\ref{pt9.1}) in $[0,T)\times \R^n$. In order to prove the lifespan estimate, we replace the initial data $(0,u_1)$ by $(0,\e f_1)$ with a small constant $\e>0$, where $f_1 \in L^1$ satisfies the assumption (\ref{optimal9.1}). Hence, there exists a sufficiently large constant $R_1 > 0$ so that we have
$$ \int_{\R^n} f_1(x)\varphi_R(x)\,dx \ge c >0 $$
for any $R > R_1$. Repeating the steps in the above proofs we arrive at the following estimate:
$$ \e \le C\, R^{-2\sigma p'+ n+ \alpha} \le C\, T^{-\frac{2\sigma p'- n- \alpha}{\alpha}} $$
with $R= T^{\frac{1}{\alpha}}$. Finally, letting $T\to T^-_\e$ we may conclude (\ref{lifetime9.1}).

\bnx
\fontshape{n}
\selectfont
We want to underline that in the special case $\sigma= 1$ an $\delta= \frac{1}{2}$ the authors in \cite{DabbiccoReissig} have investigated the critical exponent $p_{crit}=p_{crit}(n)= 1+ \frac{2}{n- 1}$. If we plug $\sigma= 1$ and $\delta= \frac{1}{2}$ into the statements of Theorem \ref{dloptimal9.1}, then the obtained results for the critical exponent $p_{crit}$ coincide.
\enx

\subsection{The case that the parameter $\sigma$ is integer and the parameter $\delta$ is fractional from $(1,\sigma)$} \label{Sec9.3}
\begin{proof}
We follow ideas from the proof of Section \ref{Sec9.2}. At first, we denote $s_\delta:= \delta- [\delta]$. Let us introduce test functions $\eta= \eta(t)$ as in Section \ref{Sec9.2} and $\varphi=\varphi(x):=\big< x\big>^{-n-2s_\delta}$. We can repeat exactly, the estimates for $J_1$ and $J_2$ as we did in the proof of Section \ref{Sec9.2} to conclude
\begin{align}
|J_1| &\lesssim I_{R,t}^{\frac{1}{p}}\, R^{-2\alpha+ \frac{n+\alpha}{p'}}, \label{t9.3.2}\\
|J_2| &\lesssim I_R^{\frac{1}{p}}\, R^{-2\sigma+ \frac{n+\alpha}{p'}}. \label{t9.3.3}
\end{align}
Let us turn to estimate $J_3$, where $\delta$ is any fractional number in $(1,\sigma)$. In the first step, applying Parseval-Plancherel formula and H\"{o}lder's inequality lead to
$$|J_3|\le I_{R,t}^{\frac{1}{p}}\, \Big(\int_{\frac{R^\alpha}{2}}^{R^{\alpha}}\int_{\R^n} \eta^{-\frac{p'}{p}}_R(t) \big|\partial_t \eta_R(t)\big|^{p'} \varphi^{-\frac{p'}{p}}_R(x)\, \big|(-\Delta)^{\delta}\varphi_R(x)\big|^{p'} \, dxdt\Big)^{\frac{1}{p'}}. $$
Now we can re-write $\delta= m_\delta+ s_\delta$, where $m_\delta:= [\delta] \ge 1$ is integer and $s_\delta$ is a fractional number in $(0,1)$. Employing Lemma \ref{bd9.3.1} we derive
$$ (-\Delta)^{\delta}\varphi_R(x)= (-\Delta)^{s_\delta} \big((-\Delta)^{m_\delta}\varphi_R(x)\big). $$
By the change of variables $\tilde{x}:= R^{-1}x$ we also notice that
$$ (-\Delta)^{m_\delta}\varphi_R(x)= R^{-2m_\delta}(-\Delta)^{m_\delta}(\varphi)(\tilde{x}) $$
since $m_\delta$ is an integer. Using the formula (\ref{RepresentationFormula^m}) we re-write
\small
\begin{align*}
(-\Delta)^{m_\delta}\varphi_R(x)&= (-1)^{m_\delta} R^{-2m_\delta} \prod_{j=0}^{m_\delta-1}(q+2j)\Big(\prod_{j=1}^{m_\delta}(-n+q+2j)\big< \tilde{x}\big>^{-q-2m_\delta} \\
&\hspace{3cm}- C^1_{m_\delta} \prod_{j=2}^{m_\delta}(-n+q+2j)(q+2m_\delta)\big< \tilde{x}\big>^{-q-2m_\delta-2} \\
&\hspace{3cm}+ C^2_{m_\delta} \prod_{j=3}^{m_\delta}(-n+q+2j)(q+2m_\delta)(q+2m_\delta+2)\big< \tilde{x}\big>^{-q-2m_\delta-4} \\
&\hspace{3cm}+\cdots+ (-1)^{m_\delta} \prod_{j=0}^{m_\delta-1}(q+2m_\delta+2j)\big< \tilde{x}\big>^{-q-4m_\delta}\Big),
\end{align*}
\normalsize
where $q:= n+2s_\delta$. For simplicity, we introduce the following functions:
$$ \varphi_k(x):= \big< x\big>^{-q-2m_\delta-2k}\quad \text{ and }\quad \varphi_{k,R}(x):= \varphi_k(R^{-1}x)=\big< \tilde{x}\big>^{-q-2m_\delta-2k} $$
with $k=0,\cdots,m_\delta$. As a result, by Lemma \ref{lemma2.3} we arrive at
\begin{align*}
(-\Delta)^{\delta}\varphi_R(x)&= (-1)^{m_\delta} R^{-2m_\delta} \prod_{j=0}^{m_\delta-1}(q+2j)\Big(\prod_{j=1}^{m_\delta}(-n+q+2j)\, (-\Delta)^{s_\delta}(\varphi_{0,R})(x) \\
&\hspace{3cm}- C^1_{m_\delta} \prod_{j=2}^{m_\delta}(-n+q+2j)(q+2m_\delta)\, (-\Delta)^{s_\delta}(\varphi_{1,R})(x) \\
&\hspace{3cm}+ C^2_{m_\delta} \prod_{j=3}^{m_\delta}(-n+q+2j)(q+2m_\delta)(q+2m_\delta+2)\, (-\Delta)^{s_\delta}(\varphi_{2,R})(x) \\
&\hspace{3cm}+\cdots+ (-1)^{m_\delta} \prod_{j=0}^{m_\delta-1}(q+2m_\delta+2j)\, (-\Delta)^{s_\delta}(\varphi_{m_\delta,R})(x)\Big) \\
&= (-1)^{m_\delta} R^{-2m_\delta-2s_\delta} \prod_{j=0}^{m_\delta-1}(q+2j)\Big(\prod_{j=1}^{m_\delta}(-n+q+2j)\, (-\Delta)^{s_\delta}(\varphi_0)(\tilde{x}) \\
&\hspace{3cm}- C^1_{m_\delta} \prod_{j=2}^{m_\delta}(-n+q+2j)(q+2m_\delta)\, (-\Delta)^{s_\delta}(\varphi_1)(\tilde{x}) \\
&\hspace{3cm}+ C^2_{m_\delta} \prod_{j=3}^{m_\delta}(-n+q+2j)(q+2m_\delta)(q+2m_\delta+2)\, (-\Delta)^{s_\delta}(\varphi_2)(\tilde{x}) \\
&\hspace{3cm}+\cdots+ (-1)^{m_\delta} \prod_{j=0}^{m_\delta-1}(q+2m_\delta+2j)\, (-\Delta)^{s_\delta}(\varphi_{m_\delta})(\tilde{x})\Big) \\
&= R^{-2\delta} (-\Delta)^{\delta}(\varphi)(\tilde{x}).
\end{align*}
For this reason, performing the change of variables $\tilde{t}:= R^{-\alpha}t$ we obtain
\begin{align*}
|J_3| &\lesssim I_{R,t}^{\frac{1}{p}}\, R^{-2\delta-\alpha+ \frac{n+\alpha}{p'}}\Big(\int_{\frac{1}{2}}^{1}\int_{\R^n} \eta^{-\frac{p'}{p}}(\tilde{t}) \big|\eta'(\tilde{t})\big|^{p'} \varphi^{-\frac{p'}{p}}(\tilde{x})\, \big|(-\Delta)^{\delta}(\varphi)(\tilde{x})\big|^{p'}\, d\tilde{x}d\tilde{t}\Big)^{\frac{1}{p'}} \\
&\lesssim I_{R,t}^{\frac{1}{p}}\, R^{-2\delta-\alpha+ \frac{n+\alpha}{p'}}\Big(\int_{\R^n} \varphi^{-\frac{p'}{p}}(\tilde{x})\, \big|(-\Delta)^{\delta}(\varphi)(\tilde{x})\big|^{p'}\, d\tilde{x}\Big)^{\frac{1}{p'}}.
\end{align*}
Here we used $\partial_t \eta_R(t)= R^{-\alpha}\eta'(\tilde{t})$ and the assumption (\ref{condition9.2.1}). After applying Lemma \ref{lemma2.2} with $q=n+2s_\delta$ and $\gamma=\delta$, i.e. $m=m_\delta$ and $s=s_\delta$, we may conclude
\begin{equation}
|J_3|\lesssim I_{R,t}^{\frac{1}{p}}\, R^{-2\delta-\alpha+ \frac{n+\alpha}{p'}}\Big(\int_{\R^n} \big< \tilde{x}\big>^{-n-2s_\delta}\, d\tilde{x}\Big)^{\frac{1}{p'}}\lesssim I_{R,t}^{\frac{1}{p}}\, R^{-2\delta-\alpha+ \frac{n+\alpha}{p'}}.
\label{t9.3.8}
\end{equation}
Finally, combining (\ref{t9.3.2}) to (\ref{t9.3.8}) and repeating arguments as in Section \ref{Sec9.2} we may complete the proof of Theorem \ref{dloptimal9.1}.
\end{proof}

\subsection{The case that the parameter $\sigma$ is fractional from $(1,\ity)$ and the parameter $\delta$ is integer} \label{Sec9.4}
\begin{proof}
We follow ideas from the proofs of Sections \ref{Sec9.2} and \ref{Sec9.3}. At first, we denote $s_\sigma:= \sigma- [\sigma]$. Let us introduce test functions $\eta= \eta(t)$ as in Section \ref{Sec9.2} and $\varphi=\varphi(x):=\big< x\big>^{-n-2s_\sigma}$. Then, repeating the proof of Sections \ref{Sec9.2} and \ref{Sec9.3} we may conclude what we wanted to prove.
\end{proof}

\subsection{The case that the parameter $\sigma$ is fractional from $(1,\ity)$ and the parameter $\delta$ is fractional from $(0,1)$} \label{Sec9.5}
\begin{proof}
We follow ideas from the proofs of Sections \ref{Sec9.2} and \ref{Sec9.4}. At first, we denote $s_\sigma:= \sigma- [\sigma]$. Next, we put $s^*:= \min\{s_\sigma,\,\delta\}$. It is obvious that $s^*$ is fractional from $(0,1)$. Let us introduce test functions $\eta= \eta(t)$ as in Section \ref{Sec9.2} and $\varphi=\varphi(x):=\big< x\big>^{-n-2s^*}$. Then, repeating the proof of Sections \ref{Sec9.2} and \ref{Sec9.4} we may conclude what we wanted to prove.
\end{proof}

\subsection{The case that the parameter $\sigma$ is fractional from $(1,\ity)$ and the parameter $\delta$ is fractional from $(1,\sigma)$} \label{Sec9.6}
\begin{proof}
We follow ideas from the proofs of Sections \ref{Sec9.2} and \ref{Sec9.5}. At first, we denote $s_\sigma:= \sigma- [\sigma]$ and $s_\delta:= \delta- [\delta]$. Next, we put $s^*:= \min\{s_\sigma,\,s_\delta\}$. It is obvious that $s^*$ is fractional from $(0,1)$. Let us introduce test functions $\eta= \eta(t)$ as in Section \ref{Sec9.2} and $\varphi=\varphi(x):=\big< x\big>^{-n-2s^*}$. Then, repeating the proof of Sections \ref{Sec9.2} and \ref{Sec9.5} we may conclude what we wanted to prove.
\end{proof}


\end{document}